\documentclass[a4paper, 11pt]{article}
\usepackage{amsmath}
\usepackage{amsfonts}
\usepackage{amssymb}
\usepackage[english]{babel}
\usepackage{graphicx}
\usepackage{amsthm,cases}
\usepackage[numbers,sort&compress]{natbib}
\usepackage[title]{appendix}

\usepackage{longtable}
\usepackage{tabularx}
\usepackage{diagbox,ulem}
\allowdisplaybreaks

\newtheorem{theorem}{Theorem}[section]

\newtheorem{lemma}[theorem]{Lemma}

\newtheorem{example}[theorem]{Example}

\newtheorem{remark}[theorem]{Remark}

\newtheorem{proposition}[theorem]{Proposition}
\newtheorem{conjecture}[theorem]{Conjecture}

\def\whitebox{{\hbox{\hskip 1pt
 \vrule height 6pt depth 1.5pt
 \lower 1.5pt\vbox to 7.5pt{\hrule width
    3.2pt\vfill\hrule width 3.2pt}%
 \vrule height 6pt depth 1.5pt
 \hskip 1pt } }}
\def\qed{\ifhmode\allowbreak\else\nobreak\fi\hfill\quad\nobreak
     \whitebox\medbreak}

\newcommand{\ignore}[1]{}
\allowdisplaybreaks

\begin{document}

\baselineskip 16pt
\title{ Non-existence of perfect binary sequences }
\author{\small  X. Niu $^1$, \ H. Cao $^1$ \thanks{Research
supported by the National Natural Science Foundation of China
under Grant 11571179, and the Priority Academic
Program Development of Jiangsu Higher Education Institutions. E-mail: caohaitao@njnu.edu.cn
} \ , and K. Feng $^2$ \thanks{Research supported by the Tsinghua National Lab. for Information Science and Technology, and NSFC with no. 11471178 and 11571007. E-mail: kfeng@math.tsinghua.edu.cn
}\\
\small $^1$ Institute of Mathematics,   Nanjing Normal University, Nanjing 210023, China\\
\small $^2$ Department of Mathematical Sciences,   Tsinghua University, Beijing 100084, China\\
}

\date{}
\maketitle

\begin{abstract}
  Binary sequences with lower autocorrelation  values
  have important applications in cryptography and communications.
  In this paper,  we present all possible parameters for binary periodical sequences
  with a 2-level autocorrelation values.
  For $n \equiv 1\pmod 4$,
  we prove some cases of Schmidt's Conjecture for perfect binary sequences. (Des. Codes Cryptogr. 78 (2016), 237-267.)
  For $n \equiv 2\pmod 4$,
  Jungnickel and Pott (Discrete Appl. Math. 95 (1999) 331-359.) left four perfect binary sequences
  as open problem and we solve three of its.
  For $n \equiv 3\pmod 4$,
  we present some nonexistence of binary sequences which all nontrivial autocorrelation values are equal 3.
  For $n \equiv 0\pmod 4$, we show that there do not exist the binary sequences
  which all nontrivial autocorrelation values are equal 4.

\bigskip

\end{abstract}

\section{Introduction}

For a binary periodical sequence ${\bf a} = (a_0,a_1, \ldots, a_{n-1}, \ldots)$ with period $n$ and $a_{j}\in \{-1,1\}$, $j\ge 0$,
 the {\it autocorrelation values} of ${\bf a}$ are defined by
$$C_{\bf a}(t) = \sum_{i=0}^{n-1}a_ ia_{i+t},\ \ \  0\le t\le n-1.$$
It is obvious that $C_{\bf a}(0) = n$, and  it is called {\it trivial} autocorrelation value.
 $C_{\bf a}(t)$, $1\le t\le n-1$, are called {\it nontrivial} or the {\it off-peak } autocorrelation values.
It is well-known that
\begin{equation}\label{tiaojian}
C_{\bf a}(t)\equiv n \pmod 4
\end{equation}
for $0\le t\le n-1$.

Binary sequences have many applications in engineering.
One of the applications is digital communication that a sequence with small aperiodic autocorrelation values
is intrinsically suited for the separation of signals from noise.
More applications details on binary sequences  may be found in \cite{B1971,G1982,TRJ1968},
and more results see \cite{ADHKM2001,B1994,CD2009,D1998,EK1992,F2017,JP1999,LS2005,LS2016,LF2016,S2016,SYF2017,SWY2017}.
In this paper, we are interest in
 binary sequences with  {\it 2-level autocorrelation values},
that is, all nontrivial autocorrelation values are equal to some constant $d$ ($C_{\bf a}(t)=d$ for $1\le t\le n-1$).
Sequences with a 2-level autocorrelation values were introduced in 1955
by Golomb who imposed this condition as one of his famous three
axioms for ``pseudo-random sequences", see \cite{G1982}.
A binary sequence with 2-level autocorrelation values is called {\it perfect}
if the nontrivial autocorrelation values $d$ are as small as possible in absolute value.

It turns out that sequences with a 2-level autocorrelation values are
equivalent to cyclic difference sets.
Let $G$ be a multiplicative abelian group of order $n$ with identity element $1_{G}$.
Let $D$ be a $k$-subset of $G$.
 The set $D$ is an $(n,k,\lambda)$ {\it difference set} ($(n,k,\lambda)$-DS) if
 every nonidentity element $g\in G$ has exactly $\lambda$ representations $g = xy^{-1}$ for
$x, y\in D$.
 If $G$ is a cyclic group,
then $D$ is an $(n,k,\lambda)$ {\it cyclic difference set} ($(n,k,\lambda)$-CDS).
  By definition, if $D$ is an $(n,k,\lambda)$-DS in $G$,
then $k(k-1)=(n-1)\lambda$
and $G\backslash D$ is an  $(n,n-k,n-2k+\lambda)$-DS in $G$.
There are more details of difference sets,  see \cite{B1971,D2014,JS1997}.

\begin{theorem}\label{dengjia0}{\rm (\cite{JP1999})}
A binary periodical sequence ${\bf a}$ with period $n$ and 2-level autocorrelation values is equivalent to an $(n,k,\lambda)$-CDS where
$C_{\bf a}(t) =d= n-4(k-\lambda), 1\le t\le n-1$.
\end{theorem}

For the perfect binary sequence,
Jungnickel and Pott \cite{JP1999} gave five different classes of cyclic difference sets corresponding to the perfect binary sequences.
The $(n,\frac{n-\sqrt{n}}{2},\frac{n-2\sqrt{n}}{4})$-CDS, $(n,\frac{n-\sqrt{2n-1}}{2},\frac{n-2\sqrt{2n-1}}{4})$-CDS,
 $(n,\frac{n-\sqrt{3n-2}}{2},\frac{n-2\sqrt{3n-2}}{4})$-CDS,
 $(n,\frac{n-\sqrt{2-n}}{2},\frac{n-2\sqrt{2-n}}{4})$-CDS and
$(n,\frac{n-1}{2},\frac{n-3}{4})$-CDS correspond to the perfect binary sequences with $d=0$, $d=1$, $d=2$, $d=-2$ and $d=-1$, respectively.
Then we give any autocorrelation value $d$ corresponding to the cyclic difference set.
The proof is similar to the Theorem \ref{dengjia0}.

\begin{lemma}\label{dengjia}
Let ${\bf a} = (a_0,a_1, \ldots, a_{n-1}, \ldots)$ be a binary periodical sequence with period $n$ and $a_{j}\in \{-1,1\}$, $j\ge 0$.
Let $G$ be a cyclic abelian group of order $n$ with $G=\langle g\rangle$ and $D=\{g^j: 0\le j\le n-1, a_j=1\}$.
Let $D^{\prime}=\{j: 0\le j\le n-1, a_j=1\}$.
 The following three cases are equivalent to each other:

\noindent 1. $C_{\bf a}(t)=d$, $1\le t\le n-1$.

\noindent 2.  $D$ is an $(n,k,\lambda)$-CDS in $G$.

\noindent 3.  $D^{\prime}$ is an $(n,k,\lambda)$-CDS in $\mathbb{Z}_n$, where $(k, \lambda)=(\frac{n+\varepsilon\sqrt{dn+n-d}}{2},
\frac{n+d+2\varepsilon\sqrt{dn+n-d}}{4})$, $\varepsilon \in \{-1, 1\}$.
\end{lemma}

By Lemma \ref{dengjia}, we have that $dn+n-d$ is a perfect square number.
Then $d \ge -1$ for $n\ge 3$ and $d = -2$ for $n= 2$ since $dn+n-d\ge 0$.
If  $n=2$ and $d=-2$, there exists the only perfect binary sequence $(-1,1,\ldots)$ \cite{JP1999}.
Then by (\ref{tiaojian}) we have

 (I) $n \equiv 1 \pmod 4$ and $d =1$.
$n=5$ and $n=13$ are the only known perfect binary sequences.
Eliahou and Kervaire \cite{EK1992},  Broughton \cite{B1994}
prove that there do not exist perfect binary sequences for $13< n< 20605$.
We prove some nonexistence of perfect binary sequences on $n \equiv 1 \pmod 4$
and give a partial solution to the conjecture posed by Schmidt \cite{S2016} in 2016.

 (II) $n \equiv 2 \pmod 4$ and $d= 2$.
$n=6$ is the only known perfect binary sequence.
Jungnickel and Pott show no perfect binary sequences for $6< n< 12545$ in \cite{JP1999}.
We prove nonexistence of three perfect binary sequences which left by Jungnickel and Pott.

 (III) $n \equiv 3 \pmod 4$ and $|d|=1$.
For $d=-1$, there are four series cyclic difference sets which construct all known perfect binary sequences.
They are Hall cyclic difference sets, Paley cyclic difference sets, twin-prime cyclic difference sets and Mersenne cyclic difference sets.
There are more details of those difference sets, see \cite{B1971,M1998}.
 We give the necessary conditions for binary sequence with $n \equiv 3 \pmod 4$  and $d=3$,
and also prove some nonexistence of them.

(IV) $n \equiv 0 \pmod 4$ and  $d \in\{0,4\}$.
If $d=0$,
there exists only one known perfect binary sequence (1, 1, 1, -1,\ldots).
It is a well known conjecture that
there do not exist perfect binary sequences with $n > 4$, in \cite{HJR1963}.
This conjecture is still open.
Leung and Schmidt
prove no optimal binary sequences for $4 < n < 548964900$ in \cite{LS2005,LS2016}.
For $d=4$, we give two binary sequences with $n=8,40$
 and also prove that a binary sequence with $n \equiv 0 \pmod 4$ and $n> 40$  does not exist.

\section{The $n \equiv 1 \pmod 4$ case}

In this section, we will present some nonexistence of perfect binary sequences  for $n \equiv 1 \pmod 4$ and $d=1$ and
also give a partial solution to the conjecture posed by Schmidt \cite{S2016} in 2016.

Perfect binary sequences are known only for $n=5$ and $n=13$, such as
$(-1,1,1,1,1,\ldots)$ and $(1,1,1,1,1,-1,-1,1,1,-1,1,-1,1,\ldots)$.
Turyn \cite{TRJ1968} reports nonexistence of those perfect binary sequences for $14\le n\le 265$.
Eliahou and Kervaire \cite{EK1992} used the result (Theorem 4.5 \cite{LES1983}) to obtained the nonexistence results for $14\le n\le 20604$, except $n=181,4901,5101,13613$.
Broughton \cite{B1994} was ruled out the four cases left by Eliahou and Kervaire \cite{EK1992}.
Then there does not exist a perfect binary sequence for $13 < n <20605$.
In 2016, Schmidt \cite{S2016} gave the following conjecture.

\begin{conjecture}\label{Conj}{\rm (\cite{S2016})}
There do not exist perfect binary sequences with $n > 13$ and $d = 1$.
\end{conjecture}

We define some concepts in number theory, and apply them to the result (Theorem 4.5 \cite{LES1983}) gives more nonexistence conclusions of the Conjecture \ref{Conj}.

Let $p$ be any prime number. For any nonzero integer $m$, let the {\it $p$-adic exponential valuation}, denoted $v_{p}(m)$,  be the highest
power of $p$ which divides $m$,  i.e.,  $v_{p}(m)=l$ if there exists a nonnegative integer $l$ such that $p^{l}|m$ and $p^{l+1}\nmid m$.
(If $m$ =0, we agree to write $v_{p}(0)=\infty$.)
Note that $v_{p}$ behaves a tittle like a logarithm would: $v_{p}(\frac{a}{b})=v_{p}(a)-v_{p}(b)$ for $a,b \in \mathbb{Z}$ and $ab\not= 0$.

Further define a map $v_{p}$ from $\mathbb{Q}$ to $\mathbb{N}\cup\{\infty\}$ as follows: For any $ \alpha,\beta\in \mathbb{Q}$,

(1) $v_{p}(\alpha) =\infty$ if $\alpha=0$.

(2) $v_{p}(\alpha\beta)=v_{p}(\alpha)+v_{p}(\beta)$. ($ \infty+ l=\infty+\infty=\infty, l\in \mathbb{N}$)

(3) $v_{p}(\alpha+\beta)\ge$ min$\{v_{p}(\alpha),v_{p}(\beta)\}$ and $v_{p}(\alpha+\beta)=$ min$\{v_{p}(\alpha),v_{p}(\beta)\}$ if $v_{p}(\alpha)\not=v_{p}(\beta)$.

\begin{proposition}\label{PEV}{\rm (\cite{H1982})}
For any nonzero integer $m\in \mathbb{Z}$ and prime number $p$,

\noindent 1. $p|m$ if and only if $v_{p}(m)\ge 1$.

\noindent 2. $p$ is a divisor of squarefree part of $m$ if and only if $v_{p}(m)$ is odd.
\end{proposition}

If $a$ and $b$ are integers, we say that $a$ is {\it semiprimitive} modulo $b$
if  there exists an  integer $c$ such that $a^c \equiv - 1 \pmod b$. Then we describe the result (Theorem 4.5 \cite{LES1983}) as follows.

\begin{theorem}\label{dj32}{\rm (\cite{LES1983})}
Suppose that there exists an $(n,k,\lambda)$-{\rm CDS}.
Let $e\ge 2$ be a divisor of $n$, and $p$ be a prime number, and
$p$ be semiprimitive modulo $e$.
Then $v_{p}(k-\lambda)$ is even.
\end{theorem}

Let $p$ be an odd prime and $a$ an integer not divisible by $p$.
Then $a$ is called a {\it quadratic residue modulo $p$} if there exists an integer $x$ such that $x^2\equiv a \pmod p$.
 We define the {\it Legendre symbol} for the odd prime $p$ as following : For any integer $a\nmid p$
 $$\biggl(\frac{a}{p}\biggr)= \begin{cases}
 ~~1 & \text{if} ~ a \text{ is a quadratic residue modulo} ~ p ,\\
 -1 & \text{if} ~ a \text{ is a quadratic nonresidue modulo} ~ p .\\
 \end{cases}$$

 \begin{theorem}\label{Lresult}{\rm (\cite{H1982})} Let $p$ be an odd prime. Then

$1. \ \bigl(\frac{2}{p}\bigr)= \begin{cases}
 ~~1 & \text{if} ~p \equiv \pm1 \pmod 8,\\
 -1 & \text{if} ~p \equiv \pm3 \pmod 8.\\
 \end{cases}$
 $\ \ 2. \ \bigl(\frac{3}{p}\bigr)= \begin{cases}
 ~~1 & \text{if} ~p \equiv \pm1 \pmod {12},\\
 -1 & \text{if}  ~p \equiv \pm5 \pmod {12}.\\
 \end{cases}$

  \end{theorem}

A perfect binary sequence with $n \equiv 1 \pmod 4$ and $d=1$ is in one-to-one correspondence with an $(n,k,\lambda)$-CDS,
where $(n,k,\lambda)=(n, \frac{1}{2}(n-\sqrt{2n-1}),\frac{1}{4}(n+1-2\sqrt{2n-1}))$.
Let $2n-1=u^2$ and $u\ge 1$.
Then $2\nmid u$ and $(n,k,\lambda)=(\frac{1}{2}(u^2+1), \frac{1}{4}(u -1)^2, \frac{1}{8}(u-1)(u-3))$.
Let $e$ be a prime divisor of $n=\frac{1}{2}(u^2+1)$.
Then we have $e$ is an odd prime and $u^2\equiv -1 \pmod e$ since $n \equiv 1 \pmod 4$.
So $-1$ is a quadratic residue modulo $e$.
Thus, we have $e \equiv 1 \pmod 4$.

In this section, we mainly apply Theorem~\ref{dj32} to get some nonexistence of perfect binary sequences.
Firstly, we consider special case $e=5$.
Since  $u^2\equiv -1 \pmod e$ and $2 \nmid u$,  we have $u\equiv 3,7 \pmod {10}$.
Then  we have  the following conclusion.

\begin{theorem}\label{T4137}
Let $u\equiv 3, 7 \pmod {10}$.
There do not exist $(\frac{1}{2}(u^2+1), \frac{1}{4}(u -1)^2, \frac{1}{8}(u-1)(u-3))$-CDS if one of the following  two conditions is valid :

\noindent 1.  $v_2(u^2-1)$ is even.

\noindent 2.  There exists a prime $p \equiv 2,3,4\pmod {5}$ such that $v_p(u+1)$ or $v_p(u-1)$ is odd.

\noindent Equivalently there do not exist perfect binary sequences with $(n,d)=(\frac{1}{2}(u^2+1),1)$.
\end{theorem}

\noindent{\bf Proof:}
If condition 1 is valid, let $p=2$ and $e=5$. Then we have $v_2(k-\lambda)=v_2(\frac{1}{8}(u^2-1))$ is odd.
Apply Theorem \ref{dj32} with $p=2$ and $e=5$ to get the conclusion.

If condition 2 is valid, let $e=5$ and $p \equiv 2,3,4\pmod {5}$.
Then $p$ is semiprimitive modulo 5.
Since $2 \nmid u$, we have gcd$(u+1,u-1)=2$.
If $v_p(u+1)$ is odd, then we have $p|(u+1)$.
So $p\nmid(u-1)$ and $v_p(u-1)$=0.
Thus, $v_2(k-\lambda)=v_2(\frac{1}{8}(u^2-1))$ is odd.
By Theorem \ref{dj32} with $e=5$, we have the conclusion.
Similarly, $v_2(k-\lambda)=v_2(\frac{1}{8}(u^2-1))$ is odd if $v_p(u-1)$ is odd.
Then we have the conclusion by Theorem \ref{dj32}.
\qed

Using elementary methods in number theory, the conditions 1 and 2 in Theorem \ref{T4137} can be expressed more explicitly.
So we have the following lemma.

\begin{lemma}\label{T41Ne}
Let $u\equiv 3, 7 \pmod {10}$ be  a positive integer. Then

\noindent 1. $v_2(u^2-1)$ is even if and only if there exists a positive integer $l\ge 1$
such that $u \equiv \pm(2^{2l+1}+(-1)^l)$  or $\pm(3\cdot2^{2l+1}+(-1)^l)\pmod {5\cdot 2^{2l+2}}$.

\noindent 2.  There exists an odd prime $p \equiv 2,3,4\pmod {5}$ such that $v_p(u+1)$ or $v_p(u-1)$ is odd if and only if
 one of the following two conditions is valid :

(2.1) $u=2bp^{2l+1} +1$, $l\ge 0$, $p\nmid b$ and
\begin{equation}\label{b-1}
b\equiv \begin{cases}
 2,~4\pmod {5} & \text{ if }  p \equiv 4\pmod {5},\\
 (-1)^{l},~(-1)^{l}\cdot2\pmod {5} & \text{ if }  p \equiv 3\pmod {5},\\
 (-1)^{l}\cdot3,~(-1)^{l}\cdot4\pmod {5} & \text{ if }  p \equiv 2\pmod {5}.\\
 \end{cases}
\end{equation}

(2.2)  $u=2bp^{2l+1} -1$, $l\ge 0$, $p\nmid b$ and
\begin{equation}\label{b-2}
b\equiv \begin{cases}
 1,~3\pmod {5} & \text{ if }  p \equiv 4\pmod {5},\\
 (-1)^{l}\cdot3,~(-1)^{l}\cdot4\pmod {5} & \text{ if }  p \equiv 3\pmod {5},\\
 (-1)^{l},~(-1)^{l}\cdot2\pmod {5} & \text{ if }  p \equiv 2\pmod {5}.\\
 \end{cases}
 \end{equation}
\end{lemma}

\noindent{\bf Proof:}
1. For any  $u\equiv 3, 7 \pmod {10}$, let $a=\frac{u-1}{2}$.
Then we have $a\equiv 1, 3 \pmod {5}$ and  $\frac{u^2-1}{4}=a(a+1)$ is a positive even number.
So $v_2(a(a+1))=v_2(a)+v_2(a+1)\ge 1$. Then,

$\hspace {0.7cm}$ $v_2(u^2-1)$ is a positive even number.

$\Longleftrightarrow$ $v_2(a)+v_2(a+1)$ is a positive even number.

$\Longleftrightarrow$ $v_2(a)$ is a positive even number or $v_2(a+1)$ is a positive even number.

Suppose  $v_2(a+1)$ is a positive even number. Then,

 $v_2(a+1)=2l$, $l\ge 1$
$\Longleftrightarrow$ $a+1=2^{2l}\cdot c$, where $2\nmid c$ $\Longleftrightarrow$ $u=2a+1=2^{2l+1}\cdot c-1$ .

\noindent Since $a\equiv 1, 3 \pmod {5}$, we have $(-1)^l\cdot c\equiv 2^{2l}\cdot c \equiv a+1 \equiv 2, 4 \pmod {5}$.
Then $c\equiv (-1)^l\cdot 7, ~ (-1)^l\cdot 9 \pmod {10}$ since $2\nmid c$.

$\hspace {0.7cm}$ $u=2^{2l+1}\cdot c-1$

$\Longleftrightarrow$ $u\equiv (-1)^l\cdot 7\cdot2^{2l+1}-1, ~ (-1)^l\cdot 9\cdot2^{2l+1}-1 \pmod {10\cdot2^{2l+1}}$

$\Longleftrightarrow$ $u\equiv (-1)^{l+1}( 3\cdot2^{2l+1}+(-1)^l), ~ (-1)^{l+1}( 2^{2l+1}+(-1)^l) \pmod {5\cdot2^{2l+2}}$.

Suppose  $v_2(a)$ is a positive even number. Then
$v_2(a)=2l$, $l\ge 1$
$\Longleftrightarrow$ $a=2^{2l}\cdot c$, $2\nmid c$.
Then $(-1)^l\cdot c\equiv 2^{2l}\cdot c \equiv a \equiv 1, 3 \pmod {5}$.
Thus, $c\equiv (-1)^l, ~ (-1)^l\cdot 3 \pmod {10}$.

$\hspace {0.7cm}$
 $u=2a+1=2^{2l+1}\cdot c+1$

$\Longleftrightarrow$ $u\equiv (-1)^{l}( 2^{2l+1}+(-1)^l), ~ (-1)^{l}( 3\cdot2^{2l+1}+(-1)^l) \pmod {5\cdot2^{2l+2}}$.

Thus, we have the conclusion.

2.  For any odd prime $p \equiv 2,3,4\pmod {5}$  and $u\ge3$, we have

$\hspace {0.7cm}$  $u\equiv 3, 7 \pmod {10}$ and $v_p(u-1)$ is odd.

$\Longleftrightarrow$ $u=c\cdot p^{2l+1} +1$ where $l\ge 0$, $p\nmid c$ and $c\cdot p^{2l+1}\equiv u-1 \equiv 2,6  \pmod {10}$.

$\Longleftrightarrow$ $u=2b\cdot p^{2l+1} +1$ where $l\ge 0$, $p\nmid b$ and $b\cdot p^{2l+1}\equiv 1,3 \pmod {5}$.

If $p \equiv 4\pmod {5}$, then we have $p^{2l+1}\equiv -1\pmod {5}$.
If $p \equiv 2,3\pmod {5}$, then we have $p^{2l+1}\equiv (-1)^{l}\cdot p\pmod {5}$.
Since $b\cdot p^{2l+1}\equiv 1,3 \pmod {5}$, we have
$$b\equiv \begin{cases}
 2,~4\pmod {5} & \text{ if }  p \equiv 4\pmod {5},\\
 (-1)^{l},~(-1)^{l}\cdot2\pmod {5} & \text{ if }  p \equiv 3\pmod {5},\\
 (-1)^{l}\cdot3,~(-1)^{l}\cdot4\pmod {5} & \text{ if }  p \equiv 2\pmod {5}.\\
 \end{cases}$$

Similarly, if $v_p(u+1)$ is odd, we have  $u=2bp^{2l+1} -1$, $l\ge 0$, $p\nmid b$ and $b$ is defined as in (\ref{b-2}).

Now, the proof is complete.\qed

Let $u \equiv \pm(2^{2l+1}+(-1)^l)$, $\pm(3\cdot2^{2l+1}+(-1)^l)\pmod {5\cdot 2^{2l+2}}$.
If $l=1$, we have $u \equiv \pm7, \pm23\pmod {80}$.
If $l=2$, we have $u \equiv \pm33, \pm97\pmod {320}$.
By Lemma \ref{T41Ne}, we have that $v_2(u^2-1)$ is even.
Applying Theorem~\ref{T4137}
we have the following example.

\begin{example}\label{T41Ex1}
Let  $u \equiv \pm7, \pm23\pmod {80}$ or  $u \equiv \pm33, \pm97\pmod {320}$.
There do not exist perfect binary sequences with $(n,d)=(\frac{u^2+1}{2},1)$.
\end{example}

 Let $u=2bp^{2l+1} +1$, $p\nmid b$ and $b$ be defined as in (\ref{b-1}).
 If $l=0$ and $p=3$, we have $u \equiv 7, 13, 43, 67\pmod {90}$.
 Let $u=2bp^{2l+1} -1$, $p\nmid b$  and $b$ be defined as in (\ref{b-2}).
 If $l=0$ and $p=3$, we have $u \equiv  23,47,77,83\pmod {90}$.
 By Lemma \ref{T41Ne}, we have that $v_3(u+1)$ or $v_3(u-1)$ is odd.
Similarly, we have the following example.

\begin{example}\label{T41Ex2}
Let  $u \equiv \pm13, \pm23, \pm43, \pm83 \pmod {90}$.
There do not exist perfect binary sequences with $(n,d)=(\frac{u^2+1}{2},1)$.
\end{example}

Next, we consider odd prime $e$ such that $e \equiv 1\pmod {4}$.
We continue to define some concepts in number theory.

For any nonzero integer $m$ and $e\nmid m$, let the {\it order of $m$ modulo $e$}, denoted $O_{e}(m)=l$,  be the least
power  of $m$ satisfying $m^{l} \equiv 1 \pmod {e}$.
It is well known that $l|\varphi(e)$ and $ \varphi(e)=e-1$, where $\varphi(e)$ is Euler function.
Let $\mathbb{F}_e$ be a finite field with $e$ elements.
Then for any element $m\in\mathbb{F}_e^*$ the order of $m$ is $O_{e}(m)$.
Let $g$ be a primitive element of $\mathbb{F}_e^*$. Then  $\mathbb{F}_e^* =\langle g\rangle$.
Let $Q_e=\langle g^2\rangle$ and $\overline{Q}_e=gQ$.
Then $Q_e$ is a set of all squared elements of $\mathbb{F}_e^*$
and $\overline{Q}_e=gQ$ is a set of all nonsquared elements of $\mathbb{F}_e^*$.

\begin{lemma}\label{T41Ne1}
Let $e$ be an odd prime satisfying $e \equiv 1\pmod {4}$ and  $m$ be an odd integer such that $|m|\ge 3$ and $m\in \overline{Q}_e$.
Then there exists an odd prime $p$ such that $v_{p}(m)$ is odd and  $p$ is semiprimitive modulo $e$.
\end{lemma}

\noindent{\bf Proof:}
Let $m=\pm p_1^{a_1}\cdot \ldots \cdot p_s^{a_s}$, where $p_1,\ldots,p_s$ are distinct primes.
Then $a_i=v_{p_i}(m)\ge 1$, $1\le i\le s$.
For any $1\le i\le s$, if $p_{i}\in Q_e$, then $p_{i}^{a_i}\in Q_e$.
If $p_{i}\in \overline{Q}_e$ and $a_i$ is even, then $p_{i}^{a_i}\in Q_e$.
Since $e \equiv 1\pmod {4}$, we have $\pm 1 \in Q_e$.
Since $m\in \overline{Q}_e$, then there exists a prime divisor $p=p_j$ of $m$ such that $p\in \overline{Q}_e$ and $a_j$ is odd.
Then $v_p(m)=a_j$ is odd.
Let $g$ be a primitive element of $\mathbb{F}_e^*$. Since $p\in \overline{Q}_e$, we have $p\equiv g^t\pmod {e}$ with $2\nmid t$.
Let $l=O_e(p)$. Then we have $1\equiv p^l \equiv g^{lt}\pmod {e}$.
So we have $(e-1)|lt$, and then $2|l$ since $2\nmid t$.
Thus, $p^{\frac{l}{2}}\equiv-1\pmod {e}$ and $p$ is semiprimitive modulo $e$.
\qed

\begin{theorem}\label{T41Ne2}
Let $e$ be a prime with $e\equiv 1\pmod {4}$. If there exists an integer $u$ satisfying the following two conditions :

\noindent 1. $2\nmid u$ and  $u^2\equiv -1\pmod {e}$.

\noindent 2.  $u\equiv 2^{l}\cdot c^2 r \pm1\pmod {2^{l+1}\cdot c^2e}$, where $c>0$, $l\ge 0$, $2|(2^l\cdot c)$, $2\nmid r$ and $r\in \overline{Q}_e$.

\noindent Then there do not exist perfect binary sequences with $(n,d)=(\frac{u^2+1}{2},1)$.
\end{theorem}

\noindent{\bf Proof:}
Since $2\nmid u$ and  $u^2\equiv -1\pmod {e}$, we have $e|n=\frac{u^2+1}{2}$.
If $u\equiv 2^{l}\cdot c^2 r +1\pmod {2^{l+1}\cdot c^2e}$, then $u-1= 2^{l}\cdot c^2 r  +2^{l+1}\cdot c^2et$ $=2^{l}\cdot c^2(r+2et), t\in \mathbb{Z}$.
Since $2\nmid r$, we have $m=r+2et$ is odd and $m\in \overline{Q}_e$.
Since $e\equiv 1\pmod {4}$ we have $\{1,-1\}\subset Q_e$. Thus, $|m|\ge 3$.
By Lemma \ref {T41Ne1}, there is an odd prime $p$ such that  $v_{p}(m)$ is odd and  $p$ is semiprimitive modulo $e$.
Then $v_p(u-1)=v_p(2^{l}\cdot c^2 m)=2v_p(c)+v_p(m)$. So $v_p(u-1)$  also is odd.
So we have $v_p(u+1)=0$ and $v_p(\frac{u^2-1}{8})=v_p(u-1)+v_p(u+1)$ is odd.
By Theorem \ref{T4137}, there do not exist $(n,k,\lambda)=(\frac{u^2+1}{2},k,\lambda)$-CDS with $k-\lambda=\frac{u^2-1}{8}$ and
there do not  exist perfect binary sequences with $(n,d)=(\frac{u^2+1}{2},1)$.

Similarly, we have the conclusion if $u\equiv 2^{l}\cdot c^2 r -1\pmod {2^{l+1}\cdot c^2e}$.
\qed

 If $e=5$, $u^2 \equiv -1\pmod {5}$ and $2\nmid u$, then we have  $u\equiv 3,7\pmod {10}$.
 Let $2\nmid r$ and $r\in \overline{Q}_e$.
 Then  $r\equiv 3,7\pmod {10}$.
Let $u\equiv 2^{l}\cdot c^2 r \pm1\pmod {2^{l+1}\cdot c^2e}$.
If  $l=c=1$, then  $u\equiv 2r \pm1 \equiv 6 \pm1, 14 \pm1  \pmod {20}$.
If $l=3$ and $c=1$, then  $u\equiv 8r \pm1 \equiv 24 \pm1, 56 \pm1  \pmod {80}$.
If $l=1$ and $c=3$, hen  $u\equiv 18r \pm1 \equiv 54 \pm1, 126 \pm1  \pmod {180}$.
Since $u\equiv 3,7\pmod {10}$, we have $u \equiv \pm7\pmod {20}$ or $u \equiv \pm23\pmod {80}$ or  $u \equiv \pm55\pmod {180}$.
Applying Theorem~\ref{T41Ne2} we have some nonexistence of perfect binary sequences.

\begin{example}\label{}
Let  $u \equiv \pm7\pmod {20}$, $u \equiv \pm23\pmod {80}$ and  $u \equiv \pm55\pmod {180}$.
Then there do not exist perfect binary sequences with $(n,d)=(\frac{u^2+1}{2},1)$.
\end{example}

 Let $e=13$, $u^2 \equiv -1\pmod {13}$ and $2\nmid u$. Then $u\equiv \pm5\pmod {26}$.
Let $l=c=1$ and $u\equiv 2^{l}\cdot c^2 r \pm1\pmod {2^{l+1}\cdot c^2e}$, where $2\nmid r$ and $r\in \overline{Q}_e$.
Then $u\equiv 2r \pm1\pmod {52}$, where  $r\equiv \pm5,\pm7,\pm11\pmod {26}$.
If  $r\equiv \pm11\pmod {26}$, then we have $u \equiv \pm21\pmod {52}$.
Applying Theorem~\ref{T41Ne2} we have some nonexistence of perfect binary sequences.

\begin{example}\label{}
Let  $u \equiv \pm21\pmod {52}$.
Then there do not exist perfect binary sequences with $(n,d)=(\frac{u^2+1}{2},1)$.
\end{example}

\section{The $n \equiv 2 \pmod 4$ case}

In this section, we shall give all possible parameters of perfect binary sequences with $n \equiv 2 \pmod 4$ and $d=2$
and also shall solve three cases of open problem in \cite{JP1999}.

Jungnickel and Pott show that  perfect binary sequences with $d=2$ do not exist for $6<n\le 12545$ \cite{JP1999},
and also gave a perfect binary sequence $(-1,1,1,1,1,1,\ldots)$.
And they left four cases $n=12546$, $n=174726$, $n=2433602$ and $n=33895686$ $(n<10^9)$ for open problem.
If we solve the four cases, we will obtain all results of $n<10^9$.

By Lemma \ref{dengjia},
a perfect binary sequence is equivalent to an $(n,k,\lambda)$-CDS,
where $(n,k,\lambda)=(n, \frac{1}{2}(n-\sqrt{3n-2}),\frac{1}{4}(n+2-2\sqrt{3n-2}))$.
So let $n=2u$ for $u\in \mathbb{N}$.
Then we have $(n,k,\lambda)=(2u, \frac{1}{2}(2u -\sqrt{6u-2}), \frac{1}{2}(u+1-\sqrt{6u-2}))$.
There are two necessary conditions for $(n,k,\lambda)$-{\rm DS},
which are known as the Bruck-Ryser-Chowla Theorem.

\begin{theorem}\label{BRC}{\rm (\cite{BR1949,CR1950})}
Suppose there exists an  $(n,k,\lambda)$-{\rm DS} in $G$ with $|G|=n$.

\noindent 1.  If $n \equiv 0 \pmod 2$,
then $k - \lambda$ is a perfect square.

\noindent 2.  If $n \equiv 1 \pmod 2$,
then there exist integers $x, y$ and $z$ (not all 0) such that
$$x^2 = (k - \lambda)y^2 + (-1)^{(n-1)/2}\lambda z^2. $$
\end{theorem}

The following result is due to  Bruck-Ryser-Chowla Theorem.

\begin{theorem}\label{DSC}{\rm (\cite{D2014})}
 If $n \equiv 1,2 \pmod 4$,
 and the square part of $n$  is divisible by a prime $p\equiv 3 \pmod 4$,
then no difference set of order $n$ exists.
\end{theorem}

In order to obtain the existence of cyclic difference sets,
we will introduce a equation and some results.
If $x^2 -dy^2 =\pm1$, we say that it is  {\it Pell equation},
and the solution of Pell equation can obtain from the quadratic filed.

\begin{theorem}\label{dj40}{\rm (\cite{H1982})}
Let $d \equiv 2,3 \pmod 4$  be a positive  integer without square divisor.
If $K=Q(\sqrt{d})$ be a real quadratic filed,
then

\noindent 1.  $\varepsilon=A +B\sqrt{d}$, $(A,B\in \mathbb{Z},A,B\ge 1)$,
where $B$ is the minimum positive integer such that $dB^2\pm 1$  is a perfect square.

\noindent 2. The all integral solutions of Pell equation
$x^2 -dy^2 =1$ are $\{(x,y)=(\pm A_k,\pm B_k): \varepsilon^k=A_k +B_k, k\ge 1\}$.
\end{theorem}

\begin{theorem}\label{dj41}{\rm (\cite{H1982})}
Let $d \equiv 1 \pmod 4$  be a positive  integer without square divisor.
If $K=Q(\sqrt{d})$ be a real quadratic filed,
then

\noindent 1.  $\varepsilon=\frac{1}{2}(A +B\sqrt{d})$, $(A,B\in \mathbb{Z},A,B\ge 1)$,
where $B$ is the minimum positive integer such that $dB^2\pm 4$  is a perfect square.

\noindent 2. The all integral solutions  of Pell equation
$x^2 -dy^2 =4$ are $\{(x,y)=(\pm A_k,\pm B_k): \varepsilon^k=A_k +B_k, k\ge 1\}$.
\end{theorem}

Next, we give all possible parameters of
$(2u, \frac{1}{2}(2u -\sqrt{6u-2}), \frac{1}{2}(u+1-\sqrt{6u-2}))$-CDS
where $u$ is a positive odd integer.

\begin{lemma}\label{021}
If there exists a $(2u, \frac{1}{2}(2u -\sqrt{6u-2}), \frac{1}{2}(u+1-\sqrt{6u-2}))$-CDS with odd integer $u\ge3$,
then $u=2B_i^2+1$,
where $\varepsilon=2+\sqrt{3}$
and $\varepsilon^i=A_i+\sqrt{3}B_i$ for $i\ge 1$.
\end{lemma}

\noindent{\bf Proof:}
Let $D$ be a $(2u, \frac{1}{2}(2u -\sqrt{6u-2}), \frac{1}{2}(u+1-\sqrt{6u-2}))$-CDS.
By Theorem \ref{BRC}, we have $k-\lambda$ is a perfect square.
So we have  $6u-2=A^2$ and
$k-\lambda=\frac{u-1}{2}=B^2$  where $A,B\ge1$.
Then we have that the equation $6u-2=A^2$ and
$ u-1=2B^2$  are equivalent to Pell equation $A^2-3B^2=1$.
By Theorem \ref{dj40},
we have $\varepsilon=2+\sqrt{3}$
and the all integral solutions of Pell equation $A^2-3B^2=1$ are $\{(A_i,B_i): \varepsilon^i=A_i+\sqrt{3}B_i, i\ge1\}$.\qed

\begin{lemma}\label{2non1}
There do not exist $(12546,6176,3040)$-CDS and
$(174726,87001,43320)$-CDS.
\end{lemma}

\noindent{\bf Proof:}
Suppose that there exists a $(12546,6176,3040)$-CDS.
Then we have $n= 2\times 3^2\times 17\times 41$ and $n \equiv 2 \pmod 4$.
Applying Theorem \ref{DSC} with $p=3$ then there does not exist a $(12546,6176,3040)$-CDS, a contradiction.

Similarly,  for $n=174726$,
 we have $n=2\times 3^2\times 17\times 571$ and $n \equiv 2 \pmod 4$.
By Theorem \ref{DSC} with $p=3$ there does not exist a $(174726,87001,43320)$-CDS.
\qed

The following result is (Corollary 1 \cite{TRJ1965}) specialised to cyclic difference sets.

\begin{theorem}\label{dj30}{\rm (\cite{TRJ1965})}
Suppose that there exists an  $(n,k,\lambda)$-{\rm CDS}.
Let $c^2$  be a divisor of $k-\lambda$ and $e\ge 2$ be a divisor of $n$.
If $e^{\prime}$ is a maximum divisor of $e$ such that $gcd(e^{\prime},c)=1$ and $c$ is semiprimitive modulo $e^{\prime}$,
then $ce\le 2^{r-1}n$, where $r$ is the number of distinct prime divisors of $gcd(e,c)$.
\end{theorem}

\begin{lemma}\label{2non2}
There does not exist a $(2433602,1215450,607050)$-CDS.
\end{lemma}

\noindent{\bf Proof:}
Suppose that there exists a $(2433602,1215450,607050)$-CDS.
Let $e=1216801$ and $c=3$.
Then we have $k-\lambda=608400$, $3^2|(k-\lambda)$ and  $e^{\prime}=1216801$ is the maximum divisor of $e$ such that gcd$(c,e^{\prime})=1$.
Thus, we have $3^{20235}\equiv -1 \pmod {1216801}$.
By Theorem \ref{dj30},
we have $c\cdot e\le 2^{1-1}n$,
i.e. $3650403\le 2433602$, a contradiction.
So there does not exist a $(2433602,1215450,607050)$-CDS.
\qed

 By Lemma~\ref{021}, for $n<10^9$, there are all possible  parameters for $(n,k,\lambda)$-CDS in the following.

\begin{center}
{\bf Table 1  Cyclic difference sets for $n \equiv 2 \pmod 4$ and $A_i, B_i$ }

\begin{tabular}{c|ccccccc}
  \hline
  $i$        & $1$    & $2$     & $3$         & $4$     & $5$      & $6$          & $7$        \\
 \hline
  $A_i$      & $2$    & $7$     & $26$        & $97$    & $362$    & $1351$       & $5042$      \\
   $B_i$     & $1$    & $4$     & $15$        & $56$    & $209$    & $780$        & $2911$     \\
 \hline
   $n$       & $6$    & $66$    & $902$       & $12546$ & $174726$  & $2433602$   & $33895686$ \\
   $k$       & $1$    & $26$    & $425$       & $6176$  & $87001$    & $1215450$   & $16942801$ \\
 $\lambda$   & $0$    & $10$    & $200$       & $3040$  & $43320$    & $607050$    & $8468880 $ \\
$k-\lambda$  & $1$    & $16$    & $225$       & $3136$  & $43681$    & $608400$    & $8473921$  \\
 \hline
$Existence$  &$\surd$ &$\times$ & $\times$    & $\times$&$\times$   & $\times$    & $?$       \\
 \hline
\end{tabular}
\end{center}

Combining Lemmas \ref{021}-\ref{2non2} and the result in \cite{JP1999},
 we obtain the following theorem.

\begin{theorem}\label{}
 Perfect binary sequences with $d=2$ do not exist for $6<n\le 33895685$.
\end{theorem}

\section{The $n \equiv 3 \pmod 4$ case}
In this section, we mainly apply Theorem~\ref{dj32} to get some nonexistence of the binary sequences with $n \equiv 3 \pmod 4$ and $d=3$.

By Lemma~\ref{dengjia},
a binary sequence with $n \equiv 3 \pmod 4$ and $d=3$ is in one-to-one correspondence with an $(n,k,\lambda)$-CDS,
where $(n,k,\lambda)=(n, \frac{1}{2}(n-\sqrt{4n-3}),\frac{1}{4}(n+3-2\sqrt{4n-3}))$.
There exists a binary sequence $(-1,1,1,1,1,1,1,\ldots)$ with $(n,d)=(7,3)$ since there exists a $(7,1,0)$-CDS.
Let  $A=\sqrt{4n-3}$. Since $n \equiv 3 \pmod 4$,
 we have $4n-3=A^2$,  $A \equiv \pm3 \pmod 8, A\ge 5$ and
 $(n,k-\lambda)=(\frac{1}{4}(A^2+3), \frac{1}{16}(A^2-9))$.

In order to apply Theorem~\ref{dj32}, let $e$ be an odd prime  divisor of $n=\frac{1}{4}(A^2+3)$.
Then $A^2\equiv -3 \pmod e$.
So we have $e=3$, $3|A$ or $\bigl(\frac{-3}{e}\bigr)=1$, $e\ge 5$.

Firstly, we  consider the case $e=3$, $A \equiv \pm3 \pmod 8$ and $3|A$.
Then we have $A\equiv 3, 21 \pmod {24}$.
For this case, we give the following conclusion.

\begin{theorem}\label{T43Ne}
Let $A\equiv \pm3 \pmod {24}$.
If there exists a prime $p \equiv 2\pmod {3}$ such that $v_p(A^2-9)$ is odd, then there does not exist a $(\frac{1}{4}(A^2+3), \frac{1}{8}(A-1)(A-3), \frac{1}{16}(A-3)(A-5))$-CDS.
 Equivalently there does not exist a binary sequence with $(n,d)=(\frac{1}{4}(A^2+3),3)$.
\end{theorem}

\noindent{\bf Proof:}
Let $e=3$ and $p\equiv 2\pmod {3}$ be a prime integer. Then $p$ is semiprimitive modulo 3 and $3|\frac{1}{4}(A^2+3)$.
If $p=2$, then we have $v_2(k-\lambda)=v_2(\frac{1}{16}(A^2-9))=v_2(A^2-9) -4$.
Since $v_p(A^2-9)$ is odd,  $v_2(k-\lambda)$ also is odd.
Apply Theorem \ref{dj32} with $p=2$ and $e=3$ to get the conclusion.
If $p \ge 5$, then $v_p(k-\lambda)=v_p(\frac{1}{16}(A^2-9))=v_p(A^2-9)$. So $v_p(k-\lambda)$ is odd.
By Theorem \ref{dj32} with $e=3$ and $p$ is semiprimitive modulo 3, we have the conclusion.
\qed

\begin{lemma}\label{L43Ne}
If $m \equiv 2 \pmod {3}$ and $m\ge 2$,
then there exists a prime $p \equiv 2 \pmod {3}$ such that $v_{p}(m)$ is odd.
\end{lemma}

\noindent{\bf Proof:}
Let $m=p_1^{a_1}\cdot \ldots \cdot p_s^{a_s}$ and $p_1, \ldots, p_s$ be distinct primes.
For $1\le i\le s$, if $p_i \equiv 1 \pmod {3}$, then  $p_i^{a_i} \equiv 1 \pmod {3}$.
If $p_i \equiv 2 \pmod {3}$ and $2|a_i$, then  $p_i^{a_i} \equiv 1 \pmod {3}$.
Since  $m \equiv 2 \pmod {3}$, we have a $j\in \{1,2,\ldots,s\}$ such that $p_j \equiv 2 \pmod {3}$ and $2\nmid a_j$. Let $p=p_j$.
Then  $v_{p}(m)=v_{p}(p^{a_j})=a_j$.
Thus, $v_{p}(m)$ is odd.
\qed

\begin{lemma}\label{L433}
Let $A \equiv 27, 45, 51, 69 \pmod {72}$.
There does not exist a binary sequences with $(n,d)=(\frac{A^2+3}{4},3)$.
\end{lemma}

\noindent{\bf Proof:}
Let $e=3$ and $n=\frac{1}{4}(A^2+3)$. Then we have $3|n$.
If $A\equiv 27 \pmod {72}$,
let $A=72l+27$, $l\ge0$.
Then $A^2-9=(72l+30)(72l+24)=9\times 16 (12l+5)(3l+1)$.
It is easy to check that $(12l+5)(3l+1)\equiv 2 \pmod {3}$.
By Lemma \ref{L43Ne}, there exists a prime $p \equiv 2 \pmod {3}$ such that $v_{p}((12l+5)(3l+1))$ is odd.
Then $v_{p}(A^2-9)$ also is odd.
 Similarly, if $A \equiv45, 51, 69 \pmod {72}$, then  $v_{p}(A^2-9)$ also is odd.
 By Theorem \ref{T43Ne}, we have the conclusion.
\qed

Secondly, we consider the case $\bigl(\frac{-3}{e}\bigr)=1$, $e\ge 5$.
Then $\bigl(\frac{-3}{e}\bigr)=1$ $\Leftrightarrow$ $\bigl(\frac{3}{e}\bigr)\bigl(\frac{-1}{e}\bigr)=1$
$\Leftrightarrow$ $\bigl(\frac{3}{e}\bigr)=(-1)^{\frac{e-1}{2}}$.
 We apply Quadratic reciprocity to $\bigl(\frac{3}{e}\bigr)=(-1)^{\frac{e-1}{2}}$ gives $\bigl(\frac{e}{3}\bigr)=1$.
 Then we have $e \equiv1\pmod{3}$.
Since $e$ is an odd prime number, we have  $e \equiv1\pmod{6}$.
So $e\ge 7$.
Let $p$ be a prime number.
If $p$ is semiprimitive modulo $e$ ($O_e(p)$ is even) and  $v_p(A^2-9)$ is odd,
by Theorem \ref{dj32}, then there does not exist a binary sequence with $(n,d)=(\frac{1}{4}(A^2+3),3)$.
So we have the following conclusion.

\begin{theorem}\label{T43Nee}
Let $e$ and $p$ be two prime numbers such that  $e \equiv1\pmod{6}$ and $p$ is semiprimitive modulo $e$.
Let $A\equiv \pm3 \pmod {8}$ such that the following  two conditions are satisfied :

\noindent 1. $A^2\equiv -3 \pmod {e}$.

\noindent 2. one of the following three conditions is valid :

(2.1) $p=2$ and  $v_2(A^2-9)$ is odd.

(2.2) $p=3$,  $v_3(A^{\prime}-1)$ or $v_3(A^{\prime}+1)$ is odd, where $A=3A^{\prime}$.

(2.3) $p\ge5$ and  $v_p(A+3)$ or $v_p(A-3)$ is odd.

\noindent Then there does not exist a binary sequence with $(n,d)=(\frac{1}{4}(A^2+3),3)$.
\end{theorem}

\noindent{\bf Proof:} Let  $(n,k-\lambda)=(\frac{1}{4}(A^2+3), \frac{1}{16}(A^2-9))$.
Since $A^2\equiv -3 \pmod {e}$ and $e \equiv1\pmod{6}$, we have $e|n$.
If $p=2$, then we have $v_2(k-\lambda)=v_2(\frac{1}{16}(A^2-9))=v_2(A^2-9) -4$.
Since $v_p(A^2-9)$ is odd,  then $v_2(k-\lambda)$ also is odd.
Apply Theorem \ref{dj32} with $p$ is semiprimitive modulo $e$ to get that there does not exist an $(n,k,\lambda)$-CDS.
So there does not exist a binary sequence with $(n,d)=(\frac{1}{4}(A^2+3),3)$ by Lemma~\ref{dengjia}.
If $p = 3$, $A=3A^{\prime}$ and $v_3(A^{\prime}-1)$ or $v_3(A^{\prime}+1)$ is odd,
then $v_3(k-\lambda)=v_3(\frac{1}{16}(A^2-9))=v_3(A-3)+v_3(A+3)= v_3(A'-1)+v_3(A'+1)+2$. So $v_3(k-\lambda)$ is odd.
Similarly, we have the conclusion.
If $p \ge 5$ and $v_p(A+3)$ or $v_p(A-3)$ is odd, then $v_p(k-\lambda)=v_p(\frac{1}{16}(A^2-9))=v_p(A^2-9)=v_p(A-3)+v_p(A+3)$. So $v_p(k-\lambda)$ is odd.
Similarly, we have the conclusion.
\qed

We obtain the following nonexistence results from Theorem~\ref{T43Nee}.

\begin{remark}\label{T43Nee1}
For the cases 1 and 2.1 in Theorem~\ref{T43Nee}, let $A=\pm3+8l$, $l\ge 1$.
Then $A^2-9=16l(4l\pm3)$ and $v_2(A^2-9)= v_2(16l(4l\pm3))= v_2(l)+v_2(4l\pm3)+4$ $=v_2(l)+4$.
 So, $v_2(A^2-9)$ is odd if and only if $v_2(l)$ is odd.
Let  $l=2^{2t+1}\cdot r$, $t\ge 0$ and $2\nmid r$.
Then  $A=\pm3 + 2^{2h}\cdot r$, where $h\ge 2$ and $2\nmid r$.
 Thus, $A^2\equiv -3 \pmod {e}$

$\hspace {1.25 cm}$ $\Longleftrightarrow$ $(2^{2h}\cdot r)^2\pm6(2^{2h}\cdot r)+12 \equiv0 \pmod {e}$

$\hspace {1.25 cm}$ $\Longleftrightarrow$ $2^{2h}\cdot r \equiv \pm3 \pm\sqrt{9-12} \pmod {e}$

$\hspace {1.25 cm}$ $\Longleftrightarrow$ $2^{2h}\cdot r \equiv \pm3 \pm s \pmod {e}$,   $s^2 \equiv -3 \pmod {e}$, $2\nmid r$.

\noindent So,  we have some $r$ mod $2e$. Thus, we can obtain $A$ from $r$ and $h$.
By Theorem \ref{Lresult}, we have $\bigl(\frac{2}{e}\bigr)=1$ if $e \equiv \pm1\pmod{8}$.
Then if $e \equiv -1\pmod{8}$, we have $2^{\frac{e-1}{2}}\equiv 1\pmod{e}$ and $\frac{e-1}{2}$ is odd.
So $2$ is not semiprimitive modulo $e$. If $e \equiv \pm1\pmod{8}$£¬we have $\bigl(\frac{2}{e}\bigr)=-1$.
Then $2$ is semiprimitive modulo $e$. Since $e \equiv1\pmod{6}$, then $2$ is semiprimitive modulo $e$ if $e \equiv 13,19\pmod{24}$, and
$2$ is not semiprimitive modulo $e$ if $e \equiv 7 \pmod{24}$.
\end{remark}

\begin{lemma}\label{L43Rmark5}
If $A\in \{61,67,77,83\}$, then there does not exist a binary sequence with $(n,d)=(\frac{A^2+3}{4},3)$.
\end{lemma}

\noindent{\bf Proof:}
Let $(h,r)\in \{(2,5), (3,1)\}$ and $A=2^{2h}\cdot r \pm3$. Then we have $A= 61,67,77,83$.
Let $e\in \{19,1123,1483,1723\}$. Then we have $19\equiv1123\equiv1483\equiv1723\pmod {24}$.
Thus, $2$ is semiprimitive modulo $e$.
Let $(n,k-\lambda)=(\frac{A^2+3}{4}, \frac{A^2-9}{16})$.
By Theorem~\ref{T43Nee} and Remark~\ref{T43Nee1},
there does not exist a binary sequence with $(n,d)=(\frac{A^2+3}{4},3)$.
\qed

\begin{remark} \label{T43Nee2}
For the cases 1 and 2.2 in Theorem~\ref{T43Nee}, we have $e \equiv1\pmod{6}$ and $A=3A^{\prime}$. Then,
$v_3(A^{\prime}-1)$ or $v_3(A^{\prime}+1)$ is odd
$\Longleftrightarrow$  $A^{\prime}= \pm1 +3^{2l+1}\cdot r$, where $l\ge 0$ and $3\nmid r$.
 Thus, $A^2\equiv -3 \pmod {e}$
 $\Longleftrightarrow$ $(3^{2l+2}\cdot r)^2\pm6(3^{2l+2}\cdot r)+12 \equiv0 \pmod {e}$

$\hspace {2.4 cm}$ $\Longleftrightarrow$ $3^{2l+2}\cdot r \equiv \pm3 \pm\sqrt{9-12} \pmod {e}$

$\hspace {2.4 cm}$ $\Longleftrightarrow$ $3^{2l+2}\cdot r \equiv \pm3 \pm s \pmod {e}$,   $s^2 \equiv -3 \pmod {e}$, $3\nmid r$.

\noindent So,  we have some $r$ mod $3e$. Thus, we can obtain $A$ from $r$ and $l$.
By Theorem \ref{Lresult},
 we have $\bigl(\frac{3}{e}\bigr)=-1$ if $e \equiv 7 \pmod{12}$. Then $3$ is semiprimitive modulo $e$.
And we also have $\bigl(\frac{3}{e}\bigr)=1$ if $e \equiv 1\pmod{12}$.
\end{remark}
Now, we shall consider a special case of Remark~\ref{T43Nee2}.

Let $e=7$ and $p=3$.
By Remark~\ref{T43Nee2}, we have $A=3A^{\prime}$ and $A^{\prime}= \pm1 +3^{2l+1}\cdot r$, $l\ge 0$, $3\nmid r$.
Since $A \equiv \pm3 \pmod {8}$ and $A^2\equiv -3 \equiv 4\pmod {7}$,
we have $A \equiv \pm5,\pm19 \pmod {56}$.
Then  $A \equiv \pm51,\pm75 \pmod {3\times56}$ since $3|A$.
If $s=2$, we have $s^2\equiv -3 \pmod {7}$. So $3^{2l+2}\cdot r \equiv \pm3 \pm 2\equiv \pm1, \pm 2 \pmod {7}$, $3\nmid r$.
If $l \equiv 0\pmod{3}$, we have $r \equiv 1,3,4,6\pmod{7}$.
Let $l=0$ and $r\equiv 8, 10 \pmod{21}$.
Then we have $A=p^{2l+2}\cdot r +3\equiv 75,93\pmod{21\times9}$.
Since $A \equiv \pm51,\pm75 \pmod {3\times56}$, we have $A\equiv 75,93\pmod{21\times9\times8}$.
Applying Theorem~\ref{T43Nee} with  $e=7$ and $p=3$
there does not exist a binary sequence with $(n,d)=(\frac{A^2+3}{4},3)$ for $A\equiv 75,93\pmod{21\times9\times8}$.
Then we have the following example.

\begin{example}\label{L43Rmark6}
There does not exist a binary sequence with $(n,d)=(\frac{A^2+3}{4},3)$ for $A \equiv 75,93\pmod{1512}$.
\end{example}

\begin{remark} For the cases 1 and 2.3 in Theorem~\ref{T43Nee}, let $p$, $e$ be  two prime numbers such that $e \equiv1\pmod{6}$, $p\ge 5$ and $p$ is semiprimitive modulo $e$.
Then,

$v_p(A+3)$ or $v_p(A-3)$ is odd
$\Longleftrightarrow$  $A= \pm3 +p^{2l+1}\cdot r$, where $l\ge 0$ and $p\nmid r$.
Thus,

$A^2\equiv -3 \pmod {e}$
 $\Longleftrightarrow$ $(p^{2l+1}\cdot r)^2\pm6(p^{2l+1}\cdot r)+12 \equiv0 \pmod {e}$

$\hspace {3.0 cm}$ $\Longleftrightarrow$ $p^{2l+1}\cdot r \equiv \pm3 \pm\sqrt{-3} \pmod {e}$

$\hspace {3.0 cm}$ $\Longleftrightarrow$ $p^{2l+1}\cdot r \equiv \pm3 \pm s \pmod {e}$,   $s^2 \equiv -3 \pmod {e}$, $p\nmid r$.

\noindent So,  we have some $r$ mod $8pe$. Thus, we can obtain $A$ from $r$ and $l$.
\end{remark}

\begin{example}\label{L43Rmark71}
Let $p=5, e=7,l=0$ and $r=8$. Then we have $A=p^{2l+1}r-3=37$.
There does not exist a binary sequence with $(n,d)=(\frac{A^2+3}{4},3)$.
\end{example}

\begin{example}\label{L43Rmark72}
Let $e=13,p_1=7, p_2=11,l=0$ and $r=8$. Then we have $A_1=p_1^{2l+1}r+3=59$ and $A_2=p_2^{2l+1}r-3=85$.
There does not exist a binary sequence with $(n,d)=(\frac{A^2+3}{4},3)$ for $A\in \{59,85\}$.
\end{example}

{\small
\begin{center}
{\bf Table 2 Cyclic difference sets for $n \equiv 3 \pmod 4$ and $A\le 100$}

\begin{tabular}{c|cccccccccccccccccccccc}
  \hline
  $A \equiv 3 \pmod 8$                & $11$ & $19$  & $27$   & $35$ & $43$   & $51$  & $59$  & $67$   & $75$  & $83$  & $91$  & $99$  \\
 \hline
   $n=\frac{1}{4}(A^2+3)$             & $31$ & $91$  & $183$  &$307$ & $463$  & $651$ & $871$ & $1123$ & $1407$ &$1723$  & $2071$ & $2451$\\
   $k=\frac{1}{8}(A-1)(A-3)$          & $10$ & $36$  & $78$   &$136$ & $210$  & $300$ & $406$ & $528$  & $666$  &$820$   & $990$ & $1176$ \\
 $\lambda=\frac{1}{16}(A-3)(A-5)$     & $3$  & $14$  & $33$   &$60$  & $95$   & $138$ & $189$ & $248$  & $315$  &$390$   & $473$  & $564$ \\
$k-\lambda=\frac{1}{16}(A^2-9)$       & $7$  & $22$  & $45$   &$76$  & $115$  & $162$ & $217$ & $280$  & $351$  &$430$   & $517$  & $612$ \\
 \hline
$Existence$               &$\times$ &$\times$ & $\times$  &$?$ &$\times$&$\times$ & $\times$&$\times$ & $\times$&$\times$ & $\times$ & $\times$ \\
 \hline
  \hline
  $A \equiv -3 \pmod 8$             & $5$ & $13$ & $21$  & $29$  & $37$  & $45$   & $53$   & $61$ & $69$   & $77$   & $85$   & $93$  \\
 \hline
   $n=\frac{1}{4}(A^2+3)$           & $7$ & $43$ & $111$ &$211$  & $343$ & $507$  &$703$   &$931$ & $1191$ & $1483$ & $1807$ & $2163$ \\
   $k=\frac{1}{8}(A-1)(A-3)$        & $1$ & $15$ & $45$  &$91$   & $153$ & $231$  &$325$   &$435$ & $561$  & $703$  & $861$  & $1035$ \\
 $\lambda=\frac{1}{16}(A-3)(A-5)$   & $0$ & $5$  & $18$  &$39$   & $68$  & $105$  &$150$   &$203$ & $264$  & $333$  & $410$  & $495$ \\
$k-\lambda=\frac{1}{16}(A^2-9)$     & $1$ & $10$ & $27$  &$52$   & $85$  & $126$  &$175$   &$232$ & $297$  & $370$  & $451$  & $540$ \\
 \hline
$Existence$         &$\surd$ &$\times$ & $\times$ &$\times$ & $\times$&$\times$ & $?$&$\times$ & $\times$&$\times$ & $\times$ & $\times$ \\
 \hline
\end{tabular}
\end{center}}

 Baumert (Table 6.1, \cite{B1971}) give the existence of  difference sets for $k\le 50$.
So cyclic difference sets for $(n,k,\lambda) \in \{(31,10,3),(43,15,5),(91,36,14),(111,45,18)\}$
can be ruled out easily by using it.
By Lemma \ref{L433}, there does not exist a $(n,k,\lambda)$-CDS, where $(n,k,\lambda) =(\frac{1}{4}(A^2+3), \frac{1}{8}(A-1)(A-3), \frac{1}{16}(A-3)(A-5))$ and $A\in\{27,45,51,69,99\}$.
We also ruled out $A\in \{37,59,61,67,75,77,83,85,93\}$
by Lemma \ref{L43Rmark5} and  Examples \ref{L43Rmark6}, \ref{L43Rmark71} and \ref{L43Rmark72}.
Then we also can obtain the nonexistence results of $(n,k,\lambda)$, where $n =\frac{1}{4}(A^2+3)$ and $A\in\{29,43,91\}$
  by  corresponding  $(e,p)$ of  Theorem \ref{dj32},
  where  $(e,p) \in \{(211, 13)$, $(463,5)$, $(109,11) \}$.

\section{The $n \equiv 0 \pmod 4$ case}

In this section, we shall prove that there do not exist the binary sequences for $n \equiv 0 \pmod 4$ and $d=4$,
except  $n=8,40$, and we also give the binary sequences for $n=8, 40$ and $d=4$.

By  (\ref{tiaojian}), we have $d \equiv n \pmod 4$.
For the case $d=0$, there are lots of nonexistence results, see \cite{B1971,G1982,LS2005,LS2016,HJR1963,TRJ1965,TRJ1968}.
So we consider the next minimum case $d=4$.
By Corollary \ref{dengjia},
we have that a binary sequence is equivalent  to an $(n,k,\lambda)$-CDS,
where $(n,k,\lambda)=(n, \frac{1}{2}(n-\sqrt{5n-4}),\frac{1}{4}(n+4-2\sqrt{5n-4}))$.
Since $n \equiv 0 \pmod 4$, we may assume that $n=4u$ for $u\in \mathbb{N}$.
Then we have $(n,k,\lambda)=(4u, 2u -\sqrt{5u-1}, u+1-\sqrt{5u-1})$.

Now, we will give all possible parameters for $(4u, 2u -\sqrt{5u-1}, u+1-\sqrt{5u-1})$-CDS by using quadratic filed.

\begin{lemma}\label{041}
If there exists a $(4u, 2u -\sqrt{5u-1}, u+1-\sqrt{5u-1})$-CDS,
then $u=B_i^2+1$,
where $\varepsilon=\frac{3+\sqrt{5}}{2}$
and $\varepsilon^i=\frac{A_i+\sqrt{5}B_i}{2}$ for $i\ge1$.
\end{lemma}

\noindent{\bf Proof:}
Suppose that there exists a $(4u, 2u -\sqrt{5u-1}, u+1-\sqrt{5u-1})$-CDS.
Since $n$ is even,  we have $u-1$ is a perfect square by Theorem \ref{BRC}.
Thus, we may assume that
 $5u-1=A^2$ and $u-1=B^2$,
 where $A,B \in \mathbb{N}^+$.
Then   it is equivalent to Pell equation $A^2-5B^2=4$.
By Theorem \ref{dj41},
we have $\varepsilon=\frac{3+\sqrt{5}}{2}$
and all integral solution of Pell equation $A^2-5B^2=4$ in set $\{(A_i,B_i): \varepsilon^i=\frac{A_i+\sqrt{5}B_i}{2}, i\ge1\}$.\qed

A $(8,1,0)$-CDS is trivial,
so we have a binary sequence $(-1,1,1,1,1,1,1,1,\ldots)$ with $(n,d)=(8,4)$.
$D=\{1,2,3,5,6,9,14,15,18,20,25,27,35\}$ is a $(40,13,4)$-CDS in $\mathbb{Z}_{40}$ \cite{B1971}.
By Corollary \ref{dengjia}, we have a binary sequence with $(n,d)=(40,4)$.
For $n\ge 41$,  we shall show that binary sequences with $d=4$ do not exist.
So we only need to prove that corresponding cyclic difference sets do not exist.

\begin{theorem}\label{dj33}{\rm (\cite{LES1983})}
Suppose that there exists an $(n,k,\lambda)$-{\rm CDS} in cyclic group $G$ and  $h\ge 1$ is a proper factor of $n=|G|$.
If there exists an integer $m$ such that $m^2|(k-\lambda)$ and $m$ is semiprimitive modulo $e=\frac{n}{h}$,
then $h\ge m$.
\end{theorem}

\begin{theorem}\label{}
 There do not exist binary sequences with $n\ge 41$ and $d = 4$.
\end{theorem}

\noindent{\bf Proof:}
Suppose, for a contradiction, that there exists a perfect
binary sequence with $n \ge 41$ and $d = 4$.
Then there exist an $(n,k,\lambda)$-CDS,
where $(n,k,\lambda)=(4u, 2u -\sqrt{5u-1}, u+1-\sqrt{5u-1})$ for  $u> 10$.
By Lemma \ref{041}, we have  $u=B_{i}^2+1$,
 $\varepsilon=\frac{3+\sqrt{5}}{2}$
and $\varepsilon^i=\frac{A_i+\sqrt{5}B_i}{2}$ for $i\ge 3$.
Then we have $B_1=1$, $B_2=3$ and $B_3=8$.
It is easy check that $B_{i} >B_j$ for $i>j$.
Thus, we have $B_i>4$ for $i\ge 3$.

Given $i\ge 3$,
 we have $k-\lambda = u-1 = B_i^2$.
 Let $n=4(B_i^2+1)$, $h=4$, and $m=B_i$.
Then $m^2 = B_i^2 \equiv -1 \pmod { B_i^2+1}$.
By Theorem \ref{dj33},
we get $h=4 \ge B_i$, a contradiction.
\qed

\end{document}